\documentclass{amsart}
\usepackage{amssymb}
\usepackage{amsmath}
\usepackage{amsfonts}

\theoremstyle{plain}

\numberwithin{equation}{section}
\input{tcilatex}

\begin{document}
\title{RANDOM FIXED POINT THEOREMS FOR LOWER SEMICONTINUOUS CONDENSING
RANDOM OPERATORS}
\author{Monica Patriche}
\maketitle

\bigskip\ \ \ \ \ \ \ \ \ \ \ \ \ \ \ \ \ \ \ \ \ \ \ \ \ \ \ \ \ \ \ \ \ \ 
\textit{University of Bucharest}

\begin{center}
\textit{Faculty of Mathematics and Computer Science}

\textit{\ 14 Academiei Street}

\textit{010014 Bucharest, Romania}{\footnotesize \\[0pt]
E-mail: \textit{monica.patriche@yahoo.com}\\[0pt]
\bigskip }
\end{center}

\textbf{Abstract.}{\small \ }{\footnotesize In this paper, we study the
existence of the random fixed points for lower semicontinuous condensing
random operators defined on Banach spaces. Our results extend corresponding
ones present in literature.}



\section{\textbf{INTRODUCTION}}

Fixed point theory has been developed in the past decades as a very powerful
tool used in the majority of mathematical applications. Some of its notable
contributions have been extended and generalized to study a wide class of
problems arising in mechanics, physics, engineering sciences, economics and
equilibrium theory etc. New results concerning the existence of the
deterministic or random fixed points were obtained, for instance, in [1-4],
[6-8], [10], [12-31].

The main aim of this work is to establish random fixed point theorems for
lower semicontinuous condensing random operators defined on Banach spaces.
Our research enables us to improve some theorems obtained recently in \cite%
{fie3}.

The rest of the paper is organized as follows. In the following section,
some notational and terminological conventions are given. We also present,
for the reader's convenience, some results on continuity and measurability
of the operators. The fixed point theorems for lower semicontinuous
condensing random operators are stated in Section 3. Section 4 presents the
conclusions of our research.

\section{\textbf{PRELIMINARIES}}

Throughout this paper, we shall use the following notation:

$2^{D}$ denotes the set of all non-empty subsets of the set $D$. If $%
D\subset Y$, where $Y$ is a topological space, cl$D$ denotes the closure of $%
D$.

For the reader's convenience, we review a few basic definitions and results
from continuity and measurability of correspondences.

Let $X$ and $Y$ be non-empty sets. The graph of $T:$\smallskip $X\rightarrow
2^{Y}$ is the set Gr$(T):=\{(x,y)\in $\smallskip $X\times Y:y\in T(x)\}.$
Let \smallskip $X$, $Y$ be topological spaces and $T:$\smallskip $%
X\rightarrow 2^{Y}$ be a correspondence. $T$ is said to be \textit{lower
semicontinuous} if, for each $z\in $\smallskip $X$ and each open set $V$ in $%
Y$ with $T(x)\cap V\neq \emptyset $, there exists an open neighborhood $U$
of $x$ in \smallskip $X$ such that $T(y)\cap V\neq \emptyset $ for each $%
y\in U$. The lower sections of $T$ are defined by $T^{-1}(y):=\{x\in $%
\smallskip $X:y\in T(x)\}$ for each $y\in Y.$

Let $(X,d)$ be a metric space. We denote $B(x,r)=\{y\in E:d(y,x)<r\}.$ If $C$
is a subset of $X,$ then, we will denote $B(C,r)=\{y\in E:d(y,C)<r\},$ where 
$d(y,C)=\inf_{x\in C}d(y,x).\smallskip $

Let now $(\Omega ,\tciFourier ,\mu )$ be a complete, finite measure space,
and $Y$ be a topological space. The correspondence $T:\Omega \rightarrow
2^{Y}$ is said to have a \textit{measurable graph} if Gr$(T)\in \tciFourier
\otimes \alpha (Y)$, where $\alpha (Y)$ denotes the Borel $\sigma $-algebra
on $Y$ and $\otimes $ denotes the product $\sigma $-algebra. The
correspondence $T:\Omega \rightarrow 2^{Y}$ is said to be \textit{lower
measurable} if, for every open subset $V$ of $Y$, the set $%
T^{-1}(V)=\{\omega \in \Omega $ $:$ $T(\omega )\cap V\neq \emptyset $\} is
an element of $\tciFourier $. This notion of measurability is also called in
literature \textit{weak measurability} or just \textit{measurability}, in
comparison with strong measurability: the correspondence $T:\Omega
\rightarrow 2^{Y}$ is said to be \textit{strong measurable} if, for every
closed subset $V$ of $Y$, the set $\{\omega \in \Omega $ $:$ $T(\omega )\cap
V\neq \emptyset $\} is an element of $\tciFourier $. In the case when $Y$ is
separable, the strong measurability coincides with the lower measurability.

Recall (see Debreu \cite{deb2}, p. 359) that if $T:\Omega \rightarrow 2^{Y}$
has a measurable graph, then $T$ is lower measurable. Furthermore, if $%
T(\cdot )$ is closed valued and lower measurable, then $T:\Omega \rightarrow
2^{Y}$ has a measurable graph.

A mapping $T:\Omega \times X\rightarrow Y$ is called a \textit{random
operator} if, for each $x\in X$, the mapping $T(\cdot ,x):\Omega \rightarrow
Y$ is measurable. Similarly, a correspondence $T:\Omega \times X\rightarrow
2^{Y}$ is also called a random operator if, for each $x\in X$, $T(\cdot
,x):\Omega \rightarrow 2^{Y}$ is measurable. A measurable mapping $\xi
:\Omega \rightarrow Y$ is called a \textit{measurable selection of the
operator} $T:\Omega \rightarrow 2^{Y}$ if $\xi (\omega )\in T(\omega )$ for
each $\omega \in \Omega $. A measurable mapping $\xi :\Omega \rightarrow Y$
is called a \textit{random fixed point} of the random operator $T:\Omega
\times X\rightarrow Y$ (or $T:\Omega \times X\rightarrow 2^{Y})$ if for
every $\omega \in \Omega $%
\c{}
$\xi (\omega )=T(\omega ,\xi (\omega ))$ (or $\xi (\omega )\in T(\omega ,\xi
(\omega ))$).

We will need the following measurable selection theorem in order to prove
our results. \smallskip

\textbf{Proposition 2.1} (Kuratowski-Ryll-Nardzewski Selection Theorem \cite%
{k}). A weakly measurable correspondence with non-empty closed values from a
measurable space into a Polish space admits a measurable selector.

\section{MAIN RESULTS}

This section is meant to extend some results established in \cite{fie3}. The
main theorems obtained in this paper are Theorem 3.1 and Theorem 3.4, which
consider lower semicontinuous condensing random operators defined on Banach
spaces. New assumptions which induce the property of $\mathcal{C}$-almost
hemicompactness are formulated and used in the statements of our theorems.

Firstly, we recall the definitions of condensing and $\mathcal{C}$-almost
hemicompact correspondences.

Let $(X,d)$ be a metric space and $E$ be a \ non-empty subset of $X.$

The correspondence $T:E\rightarrow 2^{X}$ is said to be condensing (see \cite%
{tar}), if for each subset $C$ of $E$ such that $\gamma (C)>0,$ one has $%
\gamma (T(C))<\gamma (C),$ where $T(C)=\cup _{x\in C}T(x)$ and $\gamma $ is
the Kuratowski measure of noncompactness, i.e., for each bounded subset $A$
of $E,$

$\gamma (A)=$inf$\{e>0:A$ is covered by a finite number of sets of diameter $%
\leq e\}.$

If $A$ is not a bound subset of $E,$ we assign $\gamma (A)=\infty .$

$T$ is said to be countably condensing \cite{sha14} if $T(E)$ is bounded and 
$\gamma (T(C))<\gamma (C)$ for all countably bounded sets $C$ of $E$ with $%
\gamma (C)>0$.

If $\Omega $ is any non-empty set, we say that the operator $T:\Omega \times
E\rightarrow 2^{X}$ is condensing if, for each $\omega \in \Omega ,$ the
correspondence $T(\omega ,\cdot ):E\rightarrow 2^{X}$ is
condensing.\smallskip

The mapping $T:E\rightarrow X$ is said to\textit{\ }satisfy condition (A) 
\cite{sha1} if for any sequence $(x_{n}:n\in \mathbb{N})$ in $E$ and $D\in
C(E)$ such that $d(x_{n};D)\rightarrow 0$ and $d(x_{n};T(x_{n}))\rightarrow
0 $ as $n\rightarrow \infty $, there exists $x_{0}\in D$ with $x_{0}\in
T(x_{0})$. The map $T$ is called hemicompact \cite{sha1} if each sequence $%
(x_{n}:n\in \mathbb{N})$ in $E$ has a convergent subsequence whenever $%
d(x_{n};T(x_{n}))\rightarrow 0$ as $n\rightarrow \infty $. We observe that
every continuous hemicompact map satisfies condition (A). It is also known
(see \cite{sha14}) that if $(X,d)$ is a Fr\'{e}chet space, $E$ a closed
subset of $X$ and $T:E\rightarrow X$ is a countably condensing map, then $T$
is hemicompact.

Let $E$ be a subset of $X$ and $\mathcal{C}$ be a subfamily of $2^{E}.$ We
say that $\tau _{E}$ is $\sigma -$generated by $\mathcal{C}$ (see \cite{fie3}%
), if for each $x\in E,$ $\{x\}\in \mathcal{C}$ and for each non-empty open
subset $A$ of $E,$ there exists a sequence $(C_{n};n\in \mathbb{N})$ in $%
\mathcal{C}$ such that $A=\cup _{n=0}^{\infty }C_{n}.$ If $(\Omega
,\tciFourier )$ is a measurable space, a correspondence $F:\Omega
\rightarrow 2^{E}$ is measurable if $F^{-1}(C)\in \tciFourier $ for each $%
C\in \mathcal{C}.$ In particular, if $E$ is separable, $\tau _{E}$ is $%
\sigma -$generated by all closed balls of $E$ and if $E$ is separable and
locally compact, $\tau _{E}$ is $\sigma -$generated by the family of
non-empty compact subsets of $E$.

The correspondence $T:E\rightarrow 2^{X}$ is said to be $\mathcal{C}$-almost
hemicompact (see \cite{fie3}), if $\tau _{E}$ is $\sigma -$generated by $%
\mathcal{C}$ and for each sequence $(x_{n}:n\in \mathbb{N})$ in $E$ and $%
C\in \mathcal{C}$ such that $d(x_{n},C)+h_{T}(x_{n})\rightarrow 0$ as $%
n\rightarrow \infty ,$ there exists $x\in C$ such that $h_{T}(x)=0,$ where $%
h_{T}:E\rightarrow \mathbb{R}$ is the function defined by $%
h_{T}(x)=d(x,T(x)) $ for each $x\in E.$

If $\Omega $ is any non-empty set, we say that the operator $T:\Omega \times
E\rightarrow 2^{X}$ is $\mathcal{C}$-almost hemicompact if, for each $\omega
\in \Omega ,$ the correspondence $T(\omega ,\cdot ):E\rightarrow 2^{X}$ is $%
\mathcal{C}$-almost hemicompact.

The operator $T:\Omega \times E\rightarrow 2^{X}$ is lower semicontinuous
if, for each $\omega \in \Omega ,$ the correspondence $T(\omega ,\cdot
):E\rightarrow 2^{X}$ is lower semicontinuous. \smallskip

We also present the following result in \cite{fie3}, concerning lower
semicontinuous, $\mathcal{C}-$almost hemicompact\textit{\ }operators, which
will be further extended. \smallskip

\textbf{Lemma 3.1} \textit{(Corollary 3.4 in \cite{fie3}) Let }$(X,d)$ be a
metric space, $E$\textit{\ be a complete and separable subset of }$X,$%
\textit{\ }$T:\Omega \times E\rightarrow 2^{X}$\textit{\ be a random
operator and }$\mathcal{C}\subseteq 2^{E}.$ \textit{If for each }$\omega \in
\Omega ,$\textit{\ }$T(\omega ,\cdot )$\textit{\ is lower semicontinuous, }$%
\mathcal{C}-$\textit{almost hemicompact and there exists }$x_{\omega }\in E$%
\textit{\ such that }$x_{\omega }\in $cl$T(\omega ,x_{\omega }),$ \textit{%
then, cl}$T$\textit{\ has a random fixed point.\smallskip }

The following lemmata are useful in order to prove Theorem 4.1. \smallskip

\textbf{Lemma 3.2} (Lemma 3.7 in \textit{\cite{fie3}) Let }$(X,d)$\textit{\
be a metric space, and }$(x_{n};n\in \mathbb{N})$\textit{\ and }$(y_{n};n\in 
\mathbb{N})$ \textit{be two sequences in }$X$\textit{\ such that }$%
d(x_{n},y_{n})\rightarrow 0$ as $n\rightarrow \infty .$\textit{\ Then, }$%
\gamma (A)=\gamma (B),$ \textit{where} $A=\{x_{n};n\in \mathbb{N\}}$ and $%
B=\{y_{n};n\in \mathbb{N\}}$ and $\gamma $ is the Kuratowski measure of
noncompactness\textit{.\smallskip }

\textbf{Lemma 3.3} \textit{Let }$(X,d)$\textit{\ be a metric space, }$E$%
\textit{\ be a non-empty closed separable subset of }$X,$\textit{\ }$%
\mathcal{C}$\textit{\ be the family of all closed subsets of }$E$\textit{\
such that }$\tau _{E}$\textit{\ is }$\sigma -$\textit{generated by }$%
\mathcal{C}$\textit{\ and }$T:E\rightarrow 2^{X}$\textit{\ be a condensing
correspondence which satisfies the following condition: }

$x_{0}\notin T(x_{0})$\textit{\ implies the existence of a real }$r>0$%
\textit{\ such that }$x_{0}\notin B(T(x);r)\cap B(x,r)$\textit{\ for each} $%
x\in B(x_{0},r).$\textit{\ }

\textit{Then, }$T$\textit{\ is }$\mathcal{C}-$\textit{almost hemicompact.}$%
\smallskip $

\begin{proof}
Let us consider $C\in \mathcal{C}$ and a sequence $(x_{n}:n\in \mathbb{N})$
in $E,$ for which $d(x_{n},C)+d(x_{n},T(x_{n}))\rightarrow 0$ as $%
n\rightarrow \infty .$ Let $A$ be $\{x_{n}:n\in \mathbb{N\}}$ and suppose $%
\gamma (A)>0.$ Firstly, let us denote for each $n\in \mathbb{N},$ $%
r_{n}=d(x_{n},C).$ Since $d(x_{n},T(x_{n}))\rightarrow 0$ as $n\rightarrow
\infty ,$ there exists a sequence $\{y_{n}:n\in \mathbb{N\}}$ in $X$ such
that for each $n\in \mathbb{N},$ $y_{n}\in T(x_{n})$ and hence, $%
d(x_{n},y_{n})\rightarrow 0$ as $n\rightarrow \infty .$ Lemma 3.2 implies $%
\gamma (A)=\gamma (B),$ where $B=\{y_{n}:n\in \mathbb{N}\}.$ Further, we
exploit the fact that $T$ is condensing. Therefore, we obtain $\gamma (\cup
_{n\in \mathbb{N}}T(x_{n}))<\gamma (A)=\gamma (B).$ We notice that $B\subset
\cup _{n\in \mathbb{N}}T(x_{n})$ shows that the last assertion is a
contradiction, and then, $\gamma (A)=0.$ Consequently, the sequence $%
(x_{n}:n\in \mathbb{N)}$ has a convergent subsequence $(x_{n_{k}}:k\in 
\mathbb{N)}.$ Let $x_{0}\in X$ be $x_{0}=\lim_{n_{k}\rightarrow \infty
}x_{n_{k}}.$ Since $d(x_{0},C)\leq d(x_{0},x_{n_{k}})+d(x_{n_{k}},C),$ $%
d(x_{0},C)$ must be $0$ and the closedness of $C$ implies that $x_{0}\in C.$

Further we will prove that $x_{0}\in T(x_{0}).$ Let us assume, by contrary,
that $x_{0}\notin T(x_{0}).$ Then, according to the hypotheses$,$ there
exists $r>0$ such that $x_{0}\notin B(T(x);r)\cap B(x,r)$ for each $x\in
B(x_{0},r).$ The convergence of $(x_{n_{k}}:k\in \mathbb{N}\}$ to $x_{0}$
implies the existence of a natural number $N(r)\in \mathbb{N}$ such that $%
x_{n_{k}}\in B(x_{0},r)$ for each $n_{k}>N(r).$ Consequently, $x_{0}\notin
B(T(x_{n_{k}});r)\cap B(x_{n_{k}},r)$ for each $n_{k}>N(r).$ Since for each $%
n_{k}>N(r),$ $x_{0}\in B(x_{n_{k}},r),$ it follows that if $n_{k}>N(r),$ $%
x_{0}\notin B(T(x_{n_{k}});r),$ that is $d(x_{0},T(x_{n_{k}}))>r.$ This fact
contradicts $d(x_{0},T(x_{n_{k}}))\rightarrow 0$ when $n_{k}\rightarrow
\infty ,$ which is true from the hypotheses and because $x_{0}=\lim_{n_{k}%
\rightarrow \infty }x_{n_{k}}.$ This means that our assumption is false, and
it results that $x_{0}\in T(x_{0}).$ We proved that $T$ is $\mathcal{C}-$%
almost hemicompact.
\end{proof}

\textbf{Lemma 3.4} \textit{Let }$(X,d)$\textit{\ be a metric space, }$E$%
\textit{\ be a non-empty closed separable subset of }$X,$\textit{\ }$%
\mathcal{C}$\textit{\ the family of all closed subsets of }$E$\textit{\ such
that }$\tau _{E}$\textit{\ is }$\sigma -$\textit{generated by }$\mathcal{C}$%
\textit{\ and }$T:E\rightarrow 2^{X}$\textit{\ be a correspondence such that 
}$T$\textit{\ and }$T^{-1}$\textit{\ have closed values. Therefore, if }$%
x_{0}\notin T(x_{0})$\textit{, there exists a real }$r>0$\textit{\ such that 
}$x_{0}\notin B(T(x);r)\cap B(x,r)$\textit{\ for each} $x\in B(x_{0},r).$ 
\textit{If, in addition, }$T$ \textit{is condensing, then, }$T$\textit{\ is }%
$\mathcal{C}-$\textit{almost hemicompact}$.\smallskip $

\begin{proof}
Let us consider $x_{0}\in E$ such that $x_{0}\notin T(x_{0}).$ Since $%
\{x_{0}\}\cap T^{-1}(x_{0})=\emptyset $ and $X$ is a regular space, there
exists $r_{1}>0$ such that $B(x_{0},r_{1})\cap
B(T^{-1}(x_{0});r_{1}))=\emptyset ,$ and then, $B(x_{0},r_{1})\cap
T^{-1}(x_{0})=\emptyset $. Consequently, for each $x\in B(x_{0},r_{1}),$ we
have that $x\notin T^{-1}(x_{0}),$ which is equivalent with $x_{0}\notin
T(x) $ or $\{x_{0}\}\cap T(x)=\emptyset $. The closedeness of each $T(x)$
and the regularity of $X$ imply the existence of a real number $r_{2}>0$
such that $B(x_{0},r_{2})\cap T(x)=\emptyset $ for each $x\in
B(x_{0},r_{1}), $ which implies $x_{0}\notin B(T(x);r_{2})$ for each $x\in
B(x_{0},r_{1}).$ Let $r=\min \{r_{1},r_{2}\}.$ Hence, $x_{0}\notin B(T(x);r)$
for each $x\in B(x_{0},r),$ and thus, the conclusion is fulfilled. In view
of Lemma 3.3, the last assertion is true.
\end{proof}

The next result, due to Michael, is very important in the theory of
continuous selections.\smallskip

\textbf{Lemma 3.5} \textit{(Michael} \cite{mic}). \textit{Let }$X$\textit{\
be a }$T_{1}$\textit{, paracompact space. If }$Y$\textit{\ is a Banach
space, then each lower semicontinuous convex closed valued correspondence }$%
T:X\rightarrow 2^{Y}$\textit{\ admits a continuous selection.\smallskip }

\textbf{Lemma 3.6} Let $X$ denote a nonempty closed convex subset of a
Hausdorff locally convex topological vector space $E$. If $T:X\rightarrow
2^{X}$ is condensing, then there exists a nonempty compact convex subset $K$
of $X$ such that $T(x)\subset K$ for each $x\in K.\smallskip $

By using the above lemmata, we obtain the main result of our paper, that is
Theorem 3.1, which states the existence of the random fixed points for lower
semicontinuous condensing random operators defined on Banach
spaces.\smallskip

\textbf{Theorem 3.1} \textit{Let }$(\Omega ,\mathcal{F})$\textit{\ be a
measurable space, }$E$\textit{\ be a non-empty closed convex and separable
subset of a Banach space }$X$ \textit{and let }$T:\Omega \times E\rightarrow
2^{X}$\textit{\ be a lower semicontinuous and condensing random operator
with closed and convex values. Suppose that, for each }$\omega \in \Omega ,$%
\textit{\ }$(T(\omega ,\cdot ))^{-1}:X\rightarrow 2^{E}$\textit{\ is closed
valued.}

\textit{Then, }$T$\textit{\ has a random fixed point.\smallskip }

\begin{proof}
Let $\mathcal{C}$ be a family of all closed subsets of $E$ and fix $\omega
\in \Omega .$ The correspondence $T(\omega .\cdot ):E\rightarrow 2^{X}$ is
condensing and closed valued and $(T(\omega .\cdot ))^{-1}:X\rightarrow
2^{E} $ is also closed valued and therefore, by applying Lemma 3.4, we
obtain that $T$\textit{\ }is $\mathcal{C}-$almost hemicompact$.$ In order to
apply Lemma 3.1, we will prove that the set $F(\omega ):=\{x\in E:x\in
T(\omega ,x)\}\neq \emptyset .$

According to Lemma 3.6, there exists a non-empty compact convex subset $%
K(\omega )$ of $E$ such that $T(\omega ,x)\subset K(\omega )$ for each $x\in
K(\omega ).$ Lemma 3.5 implies the existence of a continuous function $%
f_{\omega }:K(\omega )\rightarrow X$ such that $f_{\omega }(x)\in T(\omega
,x)$ for each $x\in K(\omega ).$ Since $f_{\omega }(K(\omega ))\subset
K(\omega ),$ we can apply the Brouwer-Schauder fixed point theorem and we
conclude that there exists $x_{\omega }=f_{\omega }(x_{\omega })\in T(\omega
,x_{\omega }),$ or, equivalently, $F(\omega )\neq \emptyset .$

All the assumption of Lemma 3.1 are fulfilled and then, $T$\textit{\ }has a
random fixed point.
\end{proof}

\textbf{Remark 3.1 }Theorem\ 3.1 is a stronger result than Theorem 3.9 in 
\cite{fie3}. \textit{\smallskip }

The existence of the random fixed points remains valid if for each $\omega
\in \Omega ,$ $(T(\omega .\cdot ))^{-1}:X\rightarrow 2^{E}$ is lower
semicontinuous. In this case, we establish Theorem 3.2.\smallskip

\textbf{Theorem 3.2} \textit{Let }$(\Omega ,\mathcal{F})$\textit{\ be a
measurable space, Let }$(\Omega ,\mathcal{F})$\textit{\ be a measurable
space, }$E$\textit{\ be a non-empty closed convex and separable subset of a
Banach space }$X$ \textit{and let }$T:\Omega \times E\rightarrow 2^{X}$%
\textit{\ be an operator with closed values.}

\textit{Suppose that, for each }$\omega \in \Omega ,$\textit{\ }$(T(\omega
,\cdot ))^{-1}=U(\omega ,\cdot ):X\rightarrow 2^{E}$\textit{\ is lower
semicontinuous, condensing and closed convex valued, such that }$U(\cdot ,x)$
is measurable for each $x\in X$\textit{.}

\textit{Then, }$T$\textit{\ has a random fixed point.\smallskip }

\begin{proof}
According to Theorem 3.1, there exists a measurable mapping $\xi :\Omega
\rightarrow E$ such that\textit{\ }for each $\omega \in \Omega $%
\c{}
$\xi (\omega )\in (T(\omega ,\cdot ))^{-1}(\xi (\omega )),$ that is, for
each $\omega \in \Omega $%
\c{}
$\xi (\omega )\in T(\omega ,\xi (\omega )).$ Therefore, we obtained a random
fixed point for $T.$
\end{proof}

Another main result of this section is Theorem 3.4, which involves the
condensing random operators enjoying a property which will be introduced
further.

We start with defining the following condition, which we call (*) and which
is necessary to prove the existence of random fixed points. For this
purpose, we denote Fix$(T)=\{x\in E:x\in T(x)\}.\smallskip $

\textbf{Definition 3.1} Let $(X,d)$\ be a complete metric space, $E$\ be a
non-empty closed separable subset of\textit{\ }$X,$ $\mathcal{C}$ be the
family of all closed subsets of $E$ and $Z=\{z_{n}\}$ be a countable dense
subset of $E$. We say that the correspondence $T:E\rightarrow 2^{X}$\
satisfies \textit{condition\ (*)} if, for each $C\in \mathcal{C}$ with the
property that $C\cap $Fix$(T)\neq \emptyset ,$ there exists a subsequence $%
\{z_{n_{k}}\}$ of $\{z_{n}\}$ such that $d(z_{n_{k}},C)<1/k$ and $%
d(z_{n_{k}},T(z_{n_{k}}))<1/k$ for each $k\in \mathbb{N}.$\textit{\smallskip 
}

Our work will consider a simpler assumption which implies condition (*). We
refer to condition $\alpha $ defined in \cite{pat5}$.$

Let $(X,d)$\ be a metric space and $E$ be a non-empty subset of $X.$ We say
that the correspondence $T:E\rightarrow 2^{X}$\ \textit{satisfies condition\ 
}$\alpha $ (see \cite{pat5}) if $x_{0}\in T(x_{0})$ implies that for each $%
\varepsilon >0,$ there exists an open neighborhood $U_{\varepsilon }(x_{0})$
of $x_{0}$\ such that $x_{0}\in B(T(x);\varepsilon )$\ for each $x\in
U_{\varepsilon }(x_{0}).$\ If $\Omega $\ is a non-empty set, we say that the
operator $T:\Omega \times E\rightarrow 2^{X}$\ \textit{satisfies\ condition\ 
}$\alpha $\ if, for each $\omega \in \Omega ,$\ the correspondence $T(\omega
,\cdot ):E\rightarrow 2^{X}$\ satisfies condition $\alpha .\smallskip $

Next lemma shows that condition $\alpha $ is stronger than condition (*).$%
\smallskip $

\textbf{Lemma 3.7} \textit{Let }$(X,d)$\textit{\ be a complete metric space, 
}$E$\textit{\ a non-empty closed separable subset of }$X$ \textit{and }$%
\mathcal{C}$ \textit{be the family of all closed subsets of} $E.$ \textit{If
the correspondence }$T:E\rightarrow 2^{X}$\textit{\ satisfies condition }$%
\alpha $\textit{, then, }$T$\textit{\ satisfies (*).}

\begin{proof}
Let $Z=\{z_{n}\}$ be a countable dense subset of $E.$ Let $C\in \mathcal{C}$
such that $C\cap $Fix$(T)\neq \emptyset .$ Then, there exists $x_{0}\in C$
with the property that $x_{0}\in T(x_{0}).$ According to condition $\alpha $%
, for each $k\in \mathbb{N}$, there exists an open neighborhood $%
U_{k}(x_{0}) $ of $x_{0}$\ such that $x_{0}\in B(T(x),1/k)$\ for each\textit{%
\ }$x\in U_{k}(x_{0}).$\textit{\ }Then, $x_{0}\in B(T(x),1/k)$ for each $%
x\in B(x_{0},1/k)\cap U_{k}(x_{0})$ and thus, the intersection $%
B(x_{0},1/k)\cap U_{k}(x_{0})\cap B(T(x),1/k)$ is non-empty for each $x\in
B(x_{0},1/k)\cap U_{k}(x_{0}).$ Since $B(x_{0},1/k)\cap U_{k}(x_{0})$ and $%
B(T(x),1/k)$ are open sets, $B(x_{0},1/k)\cap U_{k}(x_{0})\cap
B(T(x),1/k)\neq \{x_{0}\}.$ Therefore, for each $k\in \mathbb{N}$, we can
choose $z_{n_{k}}\in B(x_{0},1/k)\cap E\cap Z,$ $z_{n_{k}}\neq x_{0}$.
Consequently, for each $k\in \mathbb{N}$, $z_{n_{k}}\in B(C$,$1/k)\cap E\cap
Z$ and $d(z_{n_{k}},T(z_{n_{k}}))<1/k,$ that is $T$ satisfies (*).
\end{proof}

We establish the following random fixed point theorem.\smallskip

\textbf{Theorem 3.3} \textit{Let }$(\Omega ,\mathcal{F})$\textit{\ be a
measurable space, }$E$\textit{\ be a non-empty closed separable subset of a
complete metric space and let }$T:\Omega \times E\rightarrow 2^{X}$\textit{\
be a }$\mathcal{C}-$\textit{almost hemicompact} \textit{random operator
which enjoys condition (*) and is closed valued. Suppose that, for each }$%
\omega \in \Omega ,$\textit{\ the set}

$F(\omega ):=\{x\in E:x\in T(\omega ,x)\}\neq \emptyset .$

\textit{Then, }$T$\textit{\ has a random fixed point.}

\begin{proof}
Let $\mathcal{C}$ be the family of all closed subsets of $E$ and $%
Z=\{z_{n}\} $ be a countable dense subset of $E.$ Let us define $F:\Omega
\rightarrow 2^{E}$ by $F(\omega )=\{x\in E:x\in T(\omega ,x)\}.$ We notice
that $F(\omega )$ is non-empty and it is also closed, since $T(\omega ,\cdot
)$ is $\mathcal{C}-$almost hemicompact.

Let us define $h_{T}:\Omega \times E\rightarrow \mathbb{R}$ by $h_{T}(\omega
,x)=d(x,T(\omega ,x)).$ The measurability of $T(\cdot ,x),$ for each $x\in E$
implies the measurability of $h_{T}(\cdot ,x),$ for each $x\in E.$

We will prove the measurability of $F.$ In order to do this, we consider $%
C\in \mathcal{C},$ and we denote $D_{n}=\{x\in E:d(x,C)<1/n\}\cap
Z=B(C,1/n)\cap Z$ and $L(C):=\tbigcap\limits_{n=1}^{\infty
}\tbigcup\limits_{x\in D_{n}}\{\omega \in \Omega :h_{T}(\omega ,x)<1/n\}.$

$L(C)$ is measurable and we will prove further that $F^{-1}(C)=L(C).$

Firstly, let us consider $\omega \in F^{-1}(C)$ and hence there exists $%
x_{0}\in C$ such that $x_{0}\in (T(\omega ,\cdot ))^{-1}(x_{0}).$

Since $T$ satisfies condition (*)$,$ for each $k\in \mathbb{N},$ there
exists $z_{n_{k}}\in B(C,1/k)\cap Z$ such that $d(z_{n_{k}},T(\omega
,z_{n_{k}}))<1/k.$ Therefore, $\omega \in L(C)$ and then, $%
F^{-1}(C)\subseteq L(C).$

For the inverse inclusion, $L(C)\subseteq F^{-1}(C),$ let us consider $%
\omega \in L(C).$ Consequently, for each $n\geq 1,$ there exists $x_{n}\in
D_{n}$ such that $h_{T}(\omega ,x_{n})<1/n$ \ and $d(x_{n},C)<1/n.$ The
property of $C$-almost hemicompactness of $T(\omega ,\cdot )$ assures the
existence of $x\in C$ such that $h_{T}(x)=0.$ Therefore, $x\in F(\omega
)\cap C$ and $\omega \in F^{-1}(C).$

We proved that for each $C\in \mathcal{C},$ $L(C)=F^{-1}(C).$ Therefore, $F$
is measurable with non-empty closed values, and according to the Kuratowski
and Ryll-Nardzewski Proposition 2.1, $F$ has a measurable selection $\xi :%
\Omega
\rightarrow E$ such that $\xi (\omega )\in T(\omega ,(\xi ,\omega ))$ for
each $\omega \in 
\Omega
$.
\end{proof}

Based on Theorem 3.3 and Lemma 3.4, we obtain the next theorem concerning
the condensing random operators which satisfy condition\textit{\ }$\alpha $.%
\textit{\smallskip }

\textbf{Theorem 3.4 } \textit{Let }$(\Omega ,\mathcal{F})$\textit{\ be a
measurable space, }$E$\textit{\ be a non-empty closed separable subset of a
complete metric space and let }$T:\Omega \times E\rightarrow 2^{X}$\textit{\
be a condensing random operator which enjoys condition }$\alpha $\textit{\
and is closed valued, such that for each }$\omega \in \Omega ,$\textit{\ }$%
(T(\omega .\cdot ))^{-1}:X\rightarrow 2^{E}$\ \textit{is closed valued.
Supposing that, for each }$\omega \in \Omega ,$\textit{\ the set }$F(\omega
):=\{x\in E:x\in T(\omega ,x)\}\neq \emptyset ,$ \textit{then, }$T$\textit{\
has a random fixed point.}

\begin{proof}
Since $T$ is condensing and for each $\omega \in \Omega ,$ $T(\omega ,\cdot
):E\rightarrow 2^{X}$ and $(T(\omega ,\cdot ))^{-1}:X\rightarrow 2^{E}$ are
closed valued, then Lemma 3.4 implies that $T$ is $\mathcal{C}$-almost
hemicompact$.$ In order to complete the proof, we apply Theorem 3.3.
\end{proof}

\textbf{Remark 3.2 }Random fixed point theoremes for multivalued countably
condensing random operators have been obtained, for instance, in \cite{aga4}.

\section{\protect\bigskip CONCLUDING\ REMARKS}

We have proven the existence of random fixed points for condensing and lower
semicontinuous random operators defined on Banach spaces. Our study has
extended on some results which had already existed in literature. It is an
open problem to prove the existence of random fixed points for new types of
operators which satisfy weak continuity assumptions.

\bigskip

\end{document}